\documentclass[11pt]{article}
\usepackage[T1]{fontenc}
\usepackage{lmodern}
\usepackage[margin=1.1in]{geometry}
\usepackage{amsmath,amssymb,amsthm,microtype}
\usepackage[hidelinks,
  pdftitle={Hypergraph containers with deferred constraints},
  pdfauthor={Tianchi Yang}]{hyperref}

\numberwithin{equation}{section}
\theoremstyle{plain}
\newtheorem{theorem}{Theorem}[section]
\newtheorem{lemma}[theorem]{Lemma}
\newtheorem{proposition}[theorem]{Proposition}
\newtheorem{corollary}[theorem]{Corollary}
\newtheorem{claim}[theorem]{Claim}
\theoremstyle{definition}
\newtheorem{definition}[theorem]{Definition}
\theoremstyle{remark}

\newenvironment{statements}
  {\begin{enumerate}%
   }
  {\end{enumerate}}

\newcommand{\rind}{r_{\mathrm{ind}}}
\newcommand{\ex}{\operatorname{ex}}
\newcommand{\Kst}{\texorpdfstring{$K_{s,t}$}{K(s,t)}}
\newcommand{\Hact}{\mathcal H_{\mathrm{act}}}
\newcommand{\Qact}{Q_{\mathrm{act}}}

\title{Hypergraph containers with deferred constraints}
\author{Tianchi Yang\thanks{\raggedright School of Mathematics, Georgia Institute of
Technology, Atlanta, GA 30332, USA. Email: \texttt{tyang439@gatech.edu}.}}
\date{}

\begin{document}
\maketitle

\begin{abstract} 
We introduce a deferral technique for the hypergraph container method: some forbidden configurations are set aside while the containers are constructed, and are used only when counting the independent sets inside each container. We apply it to two classical enumeration problems. First, for fixed \(2\le s\le t\), the number of labelled \(K_{s,t}\)-free graphs on \(n\) vertices is at most \(2^{(\frac12(t-s+1)^{1/s}+o(1))n^{2-1/s}}\). This matches the trivial lower bound \(2^{\ex(n,K_{s,t})}\) asymptotically in the exponent whenever F\"uredi's bound on \(\ex(n,K_{s,t})\) is asymptotically sharp, in particular for \(K_{2,t}\) and \(K_{3,3}\). Second, the number of Sidon subsets of \([n]\) is \(2^{(\frac12\log_2 5+o(1))\sqrt n}\). At the end of the paper, we also outline applications to other problems: induced Ramsey numbers, sets without \(k\)-term arithmetic progressions, and list colourings of linear hypergraphs. 
\end{abstract}
\section{Introduction}\label{sec:intro}

Many questions in combinatorics ask how many discrete objects avoid a
given family of forbidden configurations, and what a typical such object
looks like: graphs with no copy of a fixed graph, sets of integers with no
solution of a given linear equation, colourings with no monochromatic
edge. These objects are the independent sets of a hypergraph whose
vertices are the elements of the ground set and whose edges are the
forbidden configurations. The method of hypergraph containers, developed
independently by Balogh, Morris and Samotij~\cite{bms} and by Saxton and
Thomason~\cite{st}, describes the independent sets of a hypergraph whose
edges are not too concentrated: all of them lie in a small number of
\emph{containers}, and each container spans few edges. For graphs, the
idea goes back to the work of Kleitman and Winston~\cite{kw} on
\(C_4\)-free graphs.

Combined with a supersaturation theorem, which says that every set
somewhat larger than the extremal size spans many forbidden
configurations, a container theorem reduces many counting problems to the
corresponding extremal problems. In this way the method has given, for
example, the number of \(C_{2\ell}\)-free graphs up to a constant factor
in the exponent~\cite{ms}, short proofs of the extremal and Ramsey
theorems for sparse random graphs and sets of integers first proved by
Conlon and Gowers~\cite{cg} and by Schacht~\cite{schacht}, bounds on the
number of sets of integers avoiding linear equations~\cite{stonline}, and
lower bounds on the list chromatic number of hypergraphs~\cite{st}; see
the survey~\cite{bmsicm} for many more. More recently, Campos and
Samotij~\cite{cs} proved new versions of the container lemma, based on
occupancy estimates for random independent sets, and a conditional form
of their approach was used by Arag\~ao, Campos, Dahia, Filipe and
Marciano~\cite{acdfm} to prove an exponential upper bound for induced
Ramsey numbers.

Container arguments usually determine such counts only up to a constant
factor in the exponent. For instance, Kleitman and Winston~\cite{kw}
showed that there are at most \(2^{1.082\,n^{3/2}}\) \(C_4\)-free graphs
on \([n]\) for large \(n\), while the subgraphs of a single extremal graph
give at least \(2^{(1/2+o(1))n^{3/2}}\). The constant is lost in two
places. A container lemma produces containers small enough for an
accurate count only if no set of two vertices lies in too many edges of
the hypergraph; and inside a final container, the trivial bound
\(2^{|C|}\) ignores the forbidden configurations that the container still
spans. For \(C_4\), both losses come from pairs of vertices with many
common neighbours, since such a pair is a diagonal of many four-cycles.

In~\cite{yang}, the author showed that the number of \(C_4\)-free graphs
on \([n]\) is \(2^{(1/2+o(1))n^{3/2}}\) by treating these four-cycles
separately; the approach was inspired by joint work with Huang and Ma on
generalized four-cycles~\cite{hmy}. The container lemma is applied only to
the four-cycles whose diagonals have small codegree, and the remaining
four-cycles are used only at the end, to count the graphs inside each
container. In this paper we develop this idea, which we call
\emph{deferral}, into a general method (Theorem~\ref{thm:deferred}). The
forbidden configurations that overlap heavily are set aside while the
containers are built. Afterwards they cut each container into small
parts, from each of which a member of the counted family takes at most
one element.

As applications, we prove an upper bound for the number of
\(K_{s,t}\)-free graphs with F\"uredi's constant in the exponent
(Theorem~\ref{thm:kst}), and we determine the number of Sidon sets up to a
factor \(2^{o(\sqrt n)}\) (Theorem~\ref{thm:sidon}). Section~\ref{sec:containers}
states the container lemma we use and proves Theorem~\ref{thm:deferred}.
Sections~\ref{sec:kst} and~\ref{sec:sidon} prove Theorems~\ref{thm:kst}
and~\ref{thm:sidon}; they are independent of each other, but they follow
the same steps, and corresponding objects are denoted by the same letters.
Section~\ref{sec:further} describes further applications, to induced
Ramsey numbers, arithmetic progressions and list colourings.

All graphs are finite and simple, all logarithms are to base \(2\), and
we omit floor and ceiling signs whenever they do not affect the argument.

\subsection{\Kst-free graphs}

For a graph \(H\), let \(\ex(n,H)\) be the maximum number of edges of an
\(H\)-free graph on \(n\) vertices. All subgraphs of an extremal graph
are \(H\)-free, so there are at least \(2^{\ex(n,H)}\) labelled
\(H\)-free graphs on \([n]\). Erd\H{o}s, Kleitman and
Rothschild~\cite{ekr} showed that this lower bound is correct in the
exponent for complete graphs, and Erd\H{o}s, Frankl and
R\"odl~\cite{efr} extended this to every non-bipartite \(H\). For
bipartite \(H\) the problem is much harder, and Morris and
Saxton~\cite{ms} showed that for \(H=C_6\) the exponent exceeds
\(\ex(n,H)\) by a constant factor.

Let \(2\le s\le t\). K\H{o}v\'ari, S\'os and Tur\'an~\cite{kst54} proved
that \(\ex(n,K_{s,t})\le\frac12(t-1)^{1/s}n^{2-1/s}+O(n)\), and
F\"uredi~\cite{furedi} improved the leading constant: with
\(c_{s,t}=\frac12(t-s+1)^{1/s}\),
\[
  \ex(n,K_{s,t})\le c_{s,t}\,n^{2-1/s}+O(n^{2-2/s}).
\]
The norm graphs of Koll\'ar, R\'onyai and Szab\'o~\cite{krs}, in the
modified form of Alon, R\'onyai and Szab\'o~\cite{ars}, show that
\(\ex(n,K_{s,t})=\Theta(n^{2-1/s})\) whenever \(t>(s-1)!\). The constant
\(c_{s,t}\) is known to be attained for \(K_{2,t}\), by a construction
of F\"uredi~\cite{furedilower}, and for \(K_{3,3}\), by a construction of
Brown~\cite{brown}.

Let \(F_{s,t}(n)\) be the number of labelled \(K_{s,t}\)-free graphs on
\([n]\). Kleitman and Winston~\cite{kw} proved that
\(\log_2F_{2,2}(n)\le1.08192\,n^{3/2}\) for large \(n\), and Balogh and
Samotij~\cite{bs} proved that \(F_{s,t}(n)=2^{O(n^{2-1/s})}\) for all
fixed \(2\le s\le t\); Morris and Saxton~\cite{ms} recovered this bound
through containers, using a balanced supersaturation theorem. For
\(s=t=2\), the author~\cite{yang} determined the leading constant,
\(\log_2F_{2,2}(n)=(\frac12+o(1))n^{3/2}\). We extend the upper bound to
all \(2\le s\le t\).

\begin{theorem}\label{thm:kst}
For fixed integers \(2\le s\le t\), as \(n\to\infty\),
\[
  \log_2F_{s,t}(n)\le\bigl(c_{s,t}+o(1)\bigr)n^{2-1/s}.
\]
\end{theorem}

Together with the lower bound \(2^{\ex(n,K_{s,t})}\),
Theorem~\ref{thm:kst} determines \(F_{s,t}(n)\) in the exponent whenever
\(\ex(n,K_{s,t})=(c_{s,t}+o(1))n^{2-1/s}\). By the constructions of
F\"uredi and Brown, we obtain the following.

\begin{corollary}\label{cor:kst}
Let \(H=K_{2,t}\) with \(t\ge2\), or \(H=K_{3,3}\). The number of
labelled \(H\)-free graphs on \([n]\) is \(2^{(1+o(1))\ex(n,H)}\).
\end{corollary}

\subsection{Sidon sets}

A set \(A\) of integers is a \emph{Sidon set} if the sums \(a+b\) with
\(a,b\in A\) and \(a\le b\) are distinct. The study of such sets began
with Sidon's work on Fourier series~\cite{sidon1932,et}. Let \(s(n)\) be
the maximum size of a Sidon subset of \([n]\). Erd\H{o}s and
Tur\'an~\cite{et} proved that \(s(n)\le\sqrt n+O(n^{1/4})\), and
Singer's difference sets~\cite{singer} give the matching lower bound
\(s(n)\ge(1-o(1))\sqrt n\), as observed by Chowla and by
Erd\H{o}s~\cite{erdos44}.

Let \(F(n)\) be the number of Sidon subsets of \([n]\). Cameron and
Erd\H{o}s~\cite{ce} proved that \(F(n)/2^{s(n)}\) is unbounded and asked
for an upper bound substantially better than the trivial one. Using the
method of Kleitman and Winston, Kohayakawa, Lee, R\"odl and
Samotij~\cite{klrs} proved that
\(\log_2F(n)\le(\log_2(33e)+o(1))\sqrt n\); they also studied Sidon sets
in sparse random subsets of \([n]\). Dellamonica, Kohayakawa, Lee, R\"odl
and Samotij~\cite{bhsets} extended the order of magnitude to
\(B_h\)-sets for every fixed \(h\ge2\). In the other direction, Saxton and
Thomason~\cite{stonline} lifted a modular Sidon set of
Ruzsa~\cite{ruzsa} to show that
\(\log_2 F(n)\ge(\frac12\log_2 5+o(1))\sqrt n\). We prove the matching
upper bound.

\begin{theorem}\label{thm:sidon}
As \(n\to\infty\),
\[
  \log_2F(n)=\bigl(\tfrac12\log_25+o(1)\bigr)\sqrt n.
\]
\end{theorem}

\subsection{Outline of the proofs}\label{sec:overview}

For \(K_{s,t}\)-free graphs, a container theorem of Morris and Saxton,
combined with a deletion argument, reduces the problem to counting
subgraphs of a graph of maximum degree \(O(n^{1-1/s})\). F\"uredi's proof
of his bound considers the common neighbourhood of a set of \(s-2\)
vertices. There a \(K_{s,t}\)-free graph induces a bipartite graph between
rows and columns in which every \(s\) rows have at most \(t-s+1\) common
columns and every two columns have at most \(t-1\) common rows. We count
such bipartite graphs inside a fixed host (Proposition~\ref{prop:local}),
and combine the counts for different \((s-2)\)-sets by Shearer's
inequality, which gives the factor \(\frac12\) in \(c_{s,t}\). The
obstruction to a sharp container count comes from pairs of columns with
many common rows, and we defer the copies of \(K_{s,t-s+2}\) through such
pairs. Inside a container, each graph being counted uses at most one
column of each small part at each row, and the number of choices is
bounded by a weighted form of the K\H{o}v\'ari--S\'os--Tur\'an argument.

For Sidon sets, an iteration of the container lemma reduces the problem
to sets \(U\subseteq[n]\) of size \(O(\sqrt n)\) spanning \(O(n)\)
additive quadruples. We then defer the quadruples whose sum or difference
has many representations in \(U\). Inside a container, a Sidon set takes
at most one element from each small part, so a part of size \(a\) offers
\(1+a\) choices. In a weighted form of the Erd\H{o}s--Tur\'an argument
the same part receives weight \(\sqrt a\), and the total weight is at
most \((1+o(1))\sqrt n\). Since \(\log(1+a)/\sqrt a\) attains its maximum
over positive integers at \(a=4\), where it equals \(\frac12\log 5\), this
gives the constant in Theorem~\ref{thm:sidon}. The value \(a=4\) matches
the lower-bound construction, in which every element has four possible
representatives.

\section{Containers with deferred constraints}\label{sec:containers}

\subsection{The container lemma}\label{sec:standard}

A hypergraph \(\mathcal H\) is a family of subsets, called edges, of a
finite vertex set; it is \(r\)-uniform if all its edges have \(r\)
elements. We write \(e(\mathcal H)\) for the number of edges,
\(\mathcal I(\mathcal H)\) for the family of independent sets, and
\(\Delta_\ell(\mathcal H)\) for the maximum number of edges containing a
given set of \(\ell\) vertices. For a set \(U\) of vertices,
\(\mathcal H[U]\) consists of the edges contained in \(U\). We shall
often use that, for \(0<\varepsilon\le1\), a set of at most \(n\)
elements has at most
\begin{equation}\label{eq:binomial}
  \sum_{j\le\varepsilon n}\binom{n}j\le\Bigl(\frac e\varepsilon\Bigr)^{\varepsilon n}
\end{equation}
subsets of size at most \(\varepsilon n\), so that recording
\(o(n)\) elements costs a factor \(2^{o(n)}\).

We use the following container lemma of Campos and Samotij.

\begin{theorem}[Campos and Samotij~\cite{cs}]\label{thm:standard}
Let \(\mathcal H\) be an \(r\)-uniform hypergraph on a finite set \(V\),
where \(r\ge2\), and let \(0<p\le1/(8r^2)\). There exist a family
\(\mathcal S\) of subsets of \(V\) and maps
\(g\colon\mathcal I(\mathcal H)\to\mathcal S\) and \(f\colon\mathcal
S\to2^V\) with the following properties.
\begin{statements}
\item For every \(I\in\mathcal I(\mathcal H)\), we have
\(g(I)\subseteq I\subseteq f(g(I))\).
\item Every \(S\in\mathcal S\) has at most \(8r^2p|V|\) elements.
\item For every \(S\in\mathcal S\), there is no hypergraph
\(\mathcal H'\subseteq\mathcal H[f(S)]\) with
\(|f(S)|\Delta_\ell(\mathcal H')<p^{\ell-1}e(\mathcal H')\) for all
\(\ell\in\{2,\dots,r\}\).
\end{statements}
\end{theorem}

Thus there are at most
\begin{equation}\label{eq:fingerprints}
  |\mathcal S|\le\Bigl(\frac ep\Bigr)^{8r^2p|V|}
\end{equation}
containers \(C=f(S)\), by~\eqref{eq:binomial}, and each of them satisfies
\begin{equation}\label{eq:degree}
  e(\mathcal H')\le|C|\max_{2\le\ell\le r}p^{1-\ell}\Delta_\ell(\mathcal H')
  \qquad\text{for every }\mathcal H'\subseteq\mathcal H[C].
\end{equation}

\subsection{Deferral}\label{sec:deferred}

Let \(\mathcal H\) be a hypergraph on a set \(V\) of \(n\) elements,
whose edges are the forbidden configurations, and let
\(\mathcal F\subseteq\mathcal I(\mathcal H)\) be the family to be
counted. Theorem~\ref{thm:standard} gives
\[
  |\mathcal F|\le\sum_{C}\bigl|\{I\in\mathcal F:I\subseteq C\}\bigr|
  \le\sum_C2^{|C|}.
\]
For a sharp bound, the containers must be small, and \(2^{|C|}\) must be
replaced by an accurate count. Both steps fail when some pairs of
vertices lie in many edges: then \(\Delta_2(\mathcal H')\) is large, the
bound~\eqref{eq:degree} says little, and the iteration stops at
containers that are too large.

We therefore split the edges of \(\mathcal H\) into \emph{active} and
\emph{deferred} ones, where the active edges, forming a hypergraph
\(\Hact\), are chosen to have small codegrees. Since
\(\Hact\subseteq\mathcal H\), we have
\(\mathcal F\subseteq\mathcal I(\mathcal H)\subseteq\mathcal I(\Hact)\),
so the containers for \(\Hact\) contain every member of \(\mathcal F\),
and~\eqref{eq:degree} is now effective. The deferred edges are not
restored to the hypergraph. Inside each container we use only a
consequence of the fact that the members of \(\mathcal F\) avoid them,
which is expressed through the following structure.

\begin{definition}\label{def:relation}
Let \(C\) be a finite set and \(L\ge2\) an integer. A \emph{relation
system} on \(C\) is a family of graphs \(\Gamma_a\) \((a\in A)\) whose
vertex sets \(C_a\) are subsets of \(C\) covering \(C\). Two distinct
vertices form a \emph{short pair} if they lie at distance less than \(L\)
in at least one of the graphs \(\Gamma_a\).
\end{definition}

In the applications, \(\Gamma_a\) joins two vertices that lie together in
many deferred edges, so a set avoiding the deferred edges contains few
short pairs. In Section~\ref{sec:sidon} there is a single relation graph
on \(C\), and in Section~\ref{sec:kst} there is one for each row. The
integer \(L\) is the local scale; in both applications
\(L=(\log n)^{1/3}\).

We partition each relation graph \(\Gamma\) as follows. Fix an ordering
of its vertices. While the graph induced by the remaining vertices has a
component with at least \(L\) vertices, start a breadth-first search at
the first vertex in such a component, using the fixed order to break
ties. The first \(L\) vertices visited induce a connected subgraph.
Remove them as a \emph{large part}. When every remaining component has
fewer than \(L\) vertices, take those components as the \emph{small
parts}. This partition has the following properties.
\begin{enumerate}
\item[(P1)] Every large part has exactly \(L\) vertices.
\item[(P2)] No edge of \(\Gamma\) joins two distinct small parts.
\item[(P3)] Every part induces a connected subgraph with at most \(L\)
vertices. In particular, any two vertices of a part are at distance less
than \(L\) in \(\Gamma\).
\end{enumerate}

Given a relation system on \(C\), partition each \(\Gamma_a\) in this way.
Let \(\mathcal T(C)\) be the family of subsets of \(C\) that avoid all
large parts and meet every small part of every \(\Gamma_a\) in at most
one vertex. A member of \(\mathcal T(C)\) is determined by its
intersections with the small parts, so
\begin{equation}\label{eq:product}
  |\mathcal T(C)|\le\prod_{a\in A}\ \prod_{P\text{ small part of }\Gamma_a}
  \bigl(1+|P|\bigr),
\end{equation}
with equality if the sets \(C_a\) are disjoint. Thus a small part \(P\)
contributes a factor \(1+|P|\), rather than the factor \(2^{|P|}\) of the
trivial bound. The following theorem is the main tool of the paper.

\begin{theorem}\label{thm:deferred}
Let \(\Hact\) be an \(r\)-uniform hypergraph on \(n\) vertices, where
\(r\ge2\), and let \(\mathcal F\subseteq\mathcal I(\Hact)\).
For \(L\ge8r^2\), let \(\mathcal C\) be the family of containers given by
Theorem~\ref{thm:standard} for \(\Hact\) with \(p=1/L\). Suppose that each
\(C\in\mathcal C\) carries a relation system with
\(\sum_a|C_a|=O(|C|)\) uniformly in \(C\), and that every
\(I\in\mathcal F\) contained in \(C\) has at most \(M\) short pairs. Then
\begin{equation}\label{eq:deferred-count}
  |\mathcal F|\le
  2^{O(n\log L/L+M\log(n+1))}
  \max_{C\in\mathcal C}|\mathcal T(C)|.
\end{equation}
The constant in the exponent depends only on \(r\) and the constant in
\(O(|C|)\).
\end{theorem}

\begin{proof}
Theorem~\ref{thm:standard} puts each \(I\in\mathcal F\) in a
container \(C\). Choose one endpoint of each short pair in \(I\), and let
\(K\) be the set of chosen vertices. There are at most \((n+1)^M\)
choices for \(K\) since \(|K|\le M\). By~(P3), \(I\setminus K\)
meets every part in at most one vertex. By~(P1), there are at most
\(\sum_a|C_a|/L=O(n/L)\) large parts, each giving at most \(1+L\) choices.
These choices, \(K\), and the remaining vertices of \(I\setminus K\),
which form a member of \(\mathcal T(C)\), determine \(I\).
Finally, \eqref{eq:fingerprints} bounds the number of containers by
\((eL)^{8r^2n/L}\). Multiplying the bounds gives~\eqref{eq:deferred-count}.
\end{proof}

The exponential factor in~\eqref{eq:deferred-count} is \(2^{o(n)}\) if
\(L\to\infty\) and \(M\log(n+1)=o(n)\). Each container also
satisfies~\eqref{eq:degree} with \(p^{1-\ell}=L^{\ell-1}\). In both
applications, the main term is therefore \(|\mathcal T(C)|\) for a
single container, which we bound in terms of the sizes of the small
parts. These bounds depend on the problem, but both rely on~(P2):
vertices of distinct small parts are not joined in the relation graphs,
so the configurations they form are active, and their number is
controlled by~\eqref{eq:degree}.

\section{\Kst-free graphs}\label{sec:kst}

In this section we prove Theorem~\ref{thm:kst}. Throughout, \(2\le s\le
t\) are fixed, and we write \(\mu=t-s+1\), so that
\(c_{s,t}=\frac12\mu^{1/s}\). Asymptotic notation refers to
\(n\to\infty\); implied constants may depend on \(s\) and \(t\), and on
other parameters only where this is indicated. For a graph \(G\) and a
set \(X\) of vertices, \(N_G(X)\) is the common neighbourhood of \(X\),
and \(N_G(v)=N_G(\{v\})\). We identify spanning subgraphs with their edge
sets, and write \(\mathcal I(G)\) for the family of \(K_{s,t}\)-free
spanning subgraphs of a graph \(G\); equivalently, \(\mathcal I(G)\) is
the family of independent sets of the hypergraph on \(E(G)\) whose edges
are the edge sets of copies of \(K_{s,t}\) in \(G\). Thus
\(F_{s,t}(n)=|\mathcal I(K_n)|\).

As for \(C_4\) in~\cite{yang}, we bound \(|\mathcal I(G)|\) by
replacing \(G\) successively with simpler graphs, losing a factor
\(2^{o(n^{2-1/s})}\) at each step. In Section~\ref{sec:kst-reduction} a
container theorem of Morris and Saxton, together with a deletion
argument, reduces Theorem~\ref{thm:kst} to graphs of maximum degree
\(O(n^{1-1/s})\) (Lemma~\ref{lem:bounded-degree} and
Proposition~\ref{prop:bounded}). Following F\"uredi's argument,
Section~\ref{sec:kst-hosts} defines bipartite hosts whose rows lie in the
common neighbourhood of a set of \(s-2\) vertices, and states the main
estimate of the section, Proposition~\ref{prop:local}, on the number of
admissible subgraphs of a host. This proposition is proved in
Sections~\ref{sec:kst-containers} and~\ref{sec:kst-weights} by means of
Theorem~\ref{thm:deferred}. In Section~\ref{sec:kst-projections} we
combine these estimates over many \((s-2)\)-sets using Shearer's
inequality, which yields the factor \(\frac12\) in \(c_{s,t}\) when
\(s\ge3\). When \(s=2\) no projection is needed, and the factor
\(\frac12\) comes instead from a bound on the number of matchings.

\subsection{Reduction to graphs of bounded degree}\label{sec:kst-reduction}

We use the counting step of the K\H{o}v\'ari--S\'os--Tur\'an theorem in
the following form.

\begin{lemma}\label{lem:kst}
Let \(A_1,\dots,A_\ell\) be subsets of an \(m\)-element set, and suppose
that every \(s\)-set lies in at most \(\nu\) of them. Then
\[
  \sum_{i=1}^{\ell}|A_i|\le (s-1)\ell+\nu^{1/s}m\ell^{1-1/s}.
\]
\end{lemma}

\begin{proof}
Let \(J\) be the set of indices \(i\) with \(|A_i|\ge s\). Since
\((x-s+1)^s\le s!\binom xs\) for every integer \(x\ge s\), counting pairs
consisting of an index \(i\) and an \(s\)-subset of \(A_i\) gives
\[
  \sum_{i\in J}(|A_i|-s+1)^s\le s!\sum_{i\in J}\binom{|A_i|}{s}
  \le\nu\,s!\binom ms\le\nu m^s.
\]
By H\"older's inequality,
\(\sum_{i\in J}(|A_i|-s+1)\le|J|^{1-1/s}(\nu m^s)^{1/s}
\le\nu^{1/s}m\ell^{1-1/s}\). Since \(|A_i|\le s-1\) for \(i\notin J\),
the lemma follows from
\[
  \sum_{i=1}^{\ell}|A_i|\le(s-1)\ell+\sum_{i\in J}(|A_i|-s+1).\qedhere
\]
\end{proof}

The first container step is a theorem of Morris and Saxton.

\begin{theorem}[Morris and Saxton~\cite{ms}]\label{thm:sparse-hosts}
There exists \(A=A(s,t)>0\) such that, for all sufficiently large
integers \(n\) and \(k\) with
\[
  k\le\frac{n^{(s-1)/(s(st-1))}}{(\log n)^{1/(s-1)}},
\]
there is a family of at most \(2^{Ak^{-(s-1)}n^{2-1/s}\log k}\) graphs
on \([n]\), each with at most \(kn^{2-1/s}\) edges, such that every
\(K_{s,t}\)-free graph on \([n]\) is contained in a member of the family.
\end{theorem}

The graphs given by Theorem~\ref{thm:sparse-hosts} have average degree
\(O(n^{1-1/s})\), but they may have vertices of much larger degree. We
delete the edges at such vertices. A \(K_{s,t}\)-free graph has few edges
at a small set of vertices, by Lemma~\ref{lem:kst}, so the deleted edges
of a graph we count can be recorded at a small cost.

\begin{lemma}\label{lem:bounded-degree}
For every \(\eta>0\) there is \(\lambda>0\) such that, for all
sufficiently large \(n\),
\[
  F_{s,t}(n)\le2^{\eta n^{2-1/s}}\max\bigl\{|\mathcal I(G)|:
  G\text{ a graph on }[n]\text{ with }\Delta(G)\le\lambda n^{1-1/s}\bigr\}.
\]
\end{lemma}

\begin{proof}
Let \(A\) be the constant in Theorem~\ref{thm:sparse-hosts}. Choose an
integer \(k\) so large that \(Ak^{-(s-1)}\log k\le\eta/2\), and then
\(\delta\in(0,1)\) with \(k\delta\log(e/\delta)\le\eta/2\). For large
\(n\), Theorem~\ref{thm:sparse-hosts} gives a family \(\mathcal G\) of
at most \(2^{\eta n^{2-1/s}/2}\) graphs on \([n]\), each with at most
\(kn^{2-1/s}\) edges, such that
\[
  \mathcal I(K_n)=\bigcup_{G\in\mathcal G}\mathcal I(G).
\]

Let \(G\in\mathcal G\). For a constant \(\lambda>0\) to be chosen
below, let \(B\) be the set of vertices of degree greater than
\(\lambda n^{1-1/s}\) in \(G\), so that \(|B|\le2kn/\lambda\), and let
\(G'\) be obtained from \(G\) by deleting the edges meeting \(B\). Then
\(\Delta(G')\le\lambda n^{1-1/s}\). Let \(I\in\mathcal I(G)\). Every
\(s\)-set of vertices lies in at most \(t-1\) of the neighbourhoods
\(N_I(w)\) with \(w\in B\), so by Lemma~\ref{lem:kst} the number of edges
of \(I\) meeting \(B\) is at most
\[
  (s-1)|B|+(t-1)^{1/s}n|B|^{1-1/s}
  =O\bigl((k/\lambda)^{1-1/s}n^{2-1/s}\bigr).
\]
Choose \(\lambda\), depending on \(k\) and \(\delta\), so large that
this is at most \(\delta kn^{2-1/s}\). As \(e(G)\le kn^{2-1/s}\), the
bound~\eqref{eq:binomial} shows that \(I\setminus E(G')\) takes at most
\((e/\delta)^{\delta kn^{2-1/s}}\le2^{\eta n^{2-1/s}/2}\) values. Since
\(I\cap E(G')\in\mathcal I(G')\), it follows that
\(|\mathcal I(G)|\le2^{\eta n^{2-1/s}/2}|\mathcal I(G')|\). Summing over
\(G\in\mathcal G\) proves the lemma.
\end{proof}

The remainder of this section is devoted to the proof of the following
bound for graphs of bounded degree.

\begin{proposition}\label{prop:bounded}
Let \(\lambda>0\). For every graph \(G\) on \([n]\) with
\(\Delta(G)\le\lambda n^{1-1/s}\),
\[
  \log|\mathcal I(G)|\le(c_{s,t}+o(1))n^{2-1/s},
\]
where the \(o(1)\) term depends on \(\lambda\) but not on \(G\).
\end{proposition}

\begin{proof}[Proof of Theorem~\ref{thm:kst} assuming
Proposition~\ref{prop:bounded}]
Let \(\eta>0\), and let \(\lambda\) be given by
Lemma~\ref{lem:bounded-degree}. By Proposition~\ref{prop:bounded}, every
graph \(G\) on \([n]\) with \(\Delta(G)\le\lambda n^{1-1/s}\) satisfies
\(\log|\mathcal I(G)|\le(c_{s,t}+o(1))n^{2-1/s}\), so
Lemma~\ref{lem:bounded-degree} gives
\(\log F_{s,t}(n)\le(c_{s,t}+\eta+o(1))n^{2-1/s}\). Since \(\eta>0\) was
arbitrary, the theorem follows.
\end{proof}

\subsection{Hosts and admissible subgraphs}\label{sec:kst-hosts}

We adapt F\"uredi's common-neighbourhood argument~\cite{furedi} to
define the bipartite graphs we count.

\begin{definition}\label{def:admissible}
A \emph{host} with \(m\) rows is a bipartite graph between a set of \(m\)
\emph{rows} and a set of at most \(n\) \emph{columns}. When \(s=2\), we
also allow a graph on \([n]\) as a host; each vertex is then regarded both
as a row and as a column, and \(m=n\). A spanning subgraph \(I\) of a host
is \emph{admissible} if
\begin{enumerate}
\item[(A1)] every \(s\)-set of rows has at most \(\mu\) common neighbours
in \(I\), and
\item[(A2)] every pair of columns has at most \(t-1\) common neighbours
in \(I\).
\end{enumerate}
We write \(\mathcal A(G)\) for the family of admissible spanning
subgraphs of a host \(G\).
\end{definition}

Every subgraph of an admissible graph is admissible. If \(s=2\) and \(G\)
is a graph on \([n]\), then (A1) and (A2) both say that any two vertices
have at most \(t-1\) common neighbours, so \(\mathcal A(G)=\mathcal I(G)\).

For \(s\ge3\), let \(I\) be a \(K_{s,t}\)-free graph and let \(X\) be a
set of \(s-2\) vertices. Take disjoint labelled copies of \(N_I(X)\)
and \(V(I)\setminus X\) as rows and columns, and join a row to a column
when their underlying vertices are adjacent in \(I\). This bipartite
graph satisfies (A1) because the vertices of any \(s\) rows already
have \(X\) as common neighbours in \(I\), leaving at most
\(t-1-(s-2)=\mu\) common neighbours outside \(X\). It satisfies (A2)
because the common rows of two columns are common neighbours in \(I\)
of the \(s\)-set formed by those two column vertices and \(X\), so
there are at most \(t-1\) such rows.

The main estimate of this section is the following.

\begin{proposition}\label{prop:local}
Let \(\lambda>0\), and let \(G\) be a host with \(m\le\lambda n^{2/s}\)
rows, each of degree at most \(\lambda n^{1-1/s}\); if \(s=2\), assume
that \(G\) is a graph on \([n]\). Then
\[
  \log|\mathcal A(G)|\le
  \begin{cases}
    \mu^{1/s}mn^{1-1/s}+o(n^{1+1/s}), & s\ge3,\\
    \tfrac12\sqrt\mu\,n^{3/2}+o(n^{3/2}), & s=2,
  \end{cases}
\]
where the error terms depend on \(\lambda\) but not on \(G\).
\end{proposition}

To see what the bound means, apply Lemma~\ref{lem:kst} to the column
neighbourhoods of an admissible graph \(I\). By (A1), \(I\) has at most
\(\mu^{1/s}mn^{1-1/s}+O(n)\) edges when \(s\ge3\), and at most
\(\frac12\sqrt\mu\,n^{3/2}+O(n)\) edges when \(s=2\), since then every
edge is counted at both of its endpoints. Thus
Proposition~\ref{prop:local} says that \(\log|\mathcal A(G)|\) is at most
the K\H{o}v\'ari--S\'os--Tur\'an bound for \(e(I)\), up to the error
term. For \(s=2\) this already proves Proposition~\ref{prop:bounded}.

\begin{proof}[Proof of Proposition~\ref{prop:bounded} for \(s=2\),
assuming Proposition~\ref{prop:local}]
Regard \(G\) itself as a host, with \(m=n\). Then \(\mathcal
A(G)=\mathcal I(G)\), and Proposition~\ref{prop:local}, applied with
\(\max\{\lambda,1\}\) in place of \(\lambda\), gives
\(\log|\mathcal I(G)|\le\frac12\sqrt\mu\,n^{3/2}+o(n^{3/2})
=(c_{2,t}+o(1))n^{3/2}\).
\end{proof}

For \(s\ge3\), Proposition~\ref{prop:bounded} is deduced from
Proposition~\ref{prop:local} in Section~\ref{sec:kst-projections}.

\subsection{Containers and small parts}\label{sec:kst-containers}

In this subsection and the next we prove Proposition~\ref{prop:local}.
Fix \(\lambda>0\) and a host \(G\) as in the proposition. Then
\(e(G)\le\lambda^2n^{1+1/s}\). Throughout the proof we set
\(L=(\log n)^{1/3}\), treated as an integer. We often use that
\(2^{L^2}=n^{o(1)}\), and that \(2^{-L}\) tends to zero faster than any
negative power of \(\log n\).

Call two columns \emph{partners} in a subgraph of \(G\) if they have at
least \(2^{-L}n^{1/s}\) common rows in it. Partners cause the large
codegrees: if \(b\) and \(b'\) are partners, then two edges \(ab\) and
\(ab'\) at the same row lie together in many copies of \(K_{s,\mu+1}\),
although an admissible graph contains at most \(t-1\) common rows of
\(b\) and \(b'\) by (A2). The proof establishes the inequalities
\begin{equation}\label{eq:kst-chain}
  |\mathcal A(G)|
  \le2^{o(n^{1+1/s})}|\mathcal A(G_0)|
  \le2^{o(n^{1+1/s})}\max_{C\in\mathcal C}|\mathcal T(C)|.
\end{equation}
The first inequality is Lemma~\ref{lem:kst-columns}: after deleting a
few edges of \(G\), we obtain a graph \(G_0\) in which no column has many
partners. The second is Lemma~\ref{lem:kst-reduction}, which applies
Theorem~\ref{thm:deferred} with the copies of \(K_{s,\mu+1}\) through
pairs of partners deferred; the hypotheses of that theorem are verified
in Lemmas~\ref{lem:kst-packings} and~\ref{lem:kst-short}. The quantity
\(|\mathcal T(C)|\) is then bounded in Section~\ref{sec:kst-weights}.

For a subgraph \(G_0\subseteq G\) with the same rows and columns, let
\(D(G_0)\) be the graph on the columns in which two columns are adjacent
if they are partners in \(G_0\).

\begin{lemma}\label{lem:kst-columns}
There is a subgraph \(G_0\) of \(G\) with the same rows and columns,
depending only on \(G\), with the following properties.
\begin{statements}
\item Every column has degree at most \(n^{1/s}\log n\) in \(G_0\).
\item The graph \(D(G_0)\) has maximum degree less than \(2^{2L}\).
\item \(|\mathcal A(G)|\le2^{o(n^{1+1/s})}|\mathcal A(G_0)|\).
\end{statements}
\end{lemma}

\begin{proof}
We delete edges of \(G\) in two rounds and show that each round deletes
\(o(n^{1+1/s})\) edges of every admissible subgraph \(I\). Since
\(e(G)=O(n^{1+1/s})\), at most \(O(n/\log n)\) columns have degree
greater than \(n^{1/s}\log n\). In the first round we delete all edges
at these columns. Their neighbourhoods in \(I\) are sets of rows, and
every \(s\)-set of rows lies in at most \(\mu\) of them by (A1). By
Lemma~\ref{lem:kst} with \(\nu=\mu\) and \(m\le\lambda n^{2/s}\), the
first round deletes \(O(n^{1+1/s}(\log n)^{-(1-1/s)})=o(n^{1+1/s})\)
edges of \(I\).

In the second round we fix an ordering of the columns and repeat the
following operation, where partners are taken in the current graph. As
long as some column has at least \(2^{2L}\) partners, take the first such
column \(b\), let \(T\) consist of its first \(2^{2L}\) partners, and
delete all edges between \(T\) and the current neighbourhood of \(b\).
Each partner loses at least \(2^{-L}n^{1/s}\) edges, so each operation
deletes at least \(2^{L}n^{1/s}\) edges of \(G\), and there are
\(O(2^{-L}n)\) operations. In one operation, a row \(a\) adjacent to \(b\)
loses its \(d_a\) edges of \(I\) to \(T\), and \(d_a\le1+\binom{d_a}2\).
Summing over these rows, and using that any two columns have at most
\(t-1\) common rows in \(I\), we see that the operation deletes at most
\(n^{1/s}\log n+(t-1)\binom{2^{2L}}{2}=n^{1/s}\log n+n^{o(1)}\) edges of
\(I\). The second round therefore deletes
\(O(2^{-L}n^{1+1/s}\log n)+n^{1+o(1)}=o(n^{1+1/s})\) edges of \(I\).

Let \(G_0\) be the final graph. Properties (a) and (b) hold by
construction, and the orderings make \(G_0\) depend only on \(G\). For
\(I\in\mathcal A(G)\), the set \(I\setminus E(G_0)\) is a subset of
\(E(G)\) with \(o(n^{1+1/s})\) elements, so by~\eqref{eq:binomial} it
takes at most \(2^{o(n^{1+1/s})}\) values. Since \(I\cap E(G_0)\) is an
admissible subgraph of \(G_0\), this proves~(c).
\end{proof}

Fix \(G_0\) as in Lemma~\ref{lem:kst-columns}, and let \(D=D(G_0)\). By
Lemma~\ref{lem:kst-columns}(c), it remains to count admissible subgraphs
of \(G_0\), and for this we apply Theorem~\ref{thm:deferred}. Let
\(\mathcal H\) be the hypergraph on \(E(G_0)\) whose edges are the edge
sets of the copies of \(K_{s,\mu+1}\) in \(G_0\) with \(s\) rows and
\(\mu+1\) columns. By (A1), every admissible subgraph of \(G_0\) is
independent in \(\mathcal H\). We call a copy \emph{active} if no two of
its \(\mu+1\) columns are partners, that is, if its columns form an
independent set in \(D\), and let \(\Hact\) consist of the edge sets of
the active copies. This hypergraph is \(r\)-uniform with \(r=s(\mu+1)\).
We apply Theorem~\ref{thm:deferred} to \(\Hact\) on \(E(G_0)\), with
\(\mathcal F=\mathcal A(G_0)\). Each container is the edge set of a
spanning subgraph of \(G_0\), and we let \(\mathcal C\) be the family of
these subgraphs.

In Section~\ref{sec:kst-weights}, columns that lie in distinct small parts
will form active copies, so we need to know that a container contains few
active copies. We measure this by disjoint packings. For an \(s\)-set
\(U\) of rows and a spanning subgraph \(C\subseteq G_0\), fix a maximum
family \(\mathcal B_U(C)\) of pairwise disjoint \((\mu+1)\)-subsets of
\(N_C(U)\), each independent in \(D\), and let
\[
  B(C)=\sum_U|\mathcal B_U(C)|,
\]
where the sum runs over all \(s\)-sets \(U\) of rows. If \(C\) is
admissible, then \(B(C)=0\) by (A1). The containers are close to
admissible in the following sense.

\begin{lemma}\label{lem:kst-packings}
Every \(C\in\mathcal C\) satisfies \(B(C)=o(n^2)\), where the error term
does not depend on \(C\).
\end{lemma}

\begin{proof}
Let \(C\in\mathcal C\), and let \(\mathcal H'\subseteq\Hact[C]\)
consist of the edge sets of the copies \(U\times S\) with
\(S\in\mathcal B_U(C)\). The edge set of a copy determines its two sides
up to interchange, so \(e(\mathcal H')\ge B(C)/2\). We claim that
\[
  \Delta_2(\mathcal H')=O\bigl((2^{-L}n^{1/s})^{s-1}\bigr).
\]
Consider two edges of \(G_0\) in a copy, oriented from rows \(a_1,a_2\)
to columns \(b_1,b_2\); there are at most four orientations. If
\(a_1=a_2\), then \(b_1\) and \(b_2\) are distinct columns of a set
independent in \(D\), so they have fewer than \(2^{-L}n^{1/s}\) common
rows in \(G_0\), and the other \(s-1\) rows of the copy can be chosen in
at most \((2^{-L}n^{1/s})^{s-1}\) ways. If \(a_1\neq a_2\), the remaining
\(s-2\) rows are neighbours of \(b_1\), giving at most
\((n^{1/s}\log n)^{s-2}=o\bigl((2^{-L}n^{1/s})^{s-1}\bigr)\) choices by
Lemma~\ref{lem:kst-columns}(a). In either case the rows of the copy and
the column \(b_1\) determine at most one member of the corresponding
packing, which proves the claim.

Since \(\Delta_\ell(\mathcal H')\le\Delta_2(\mathcal H')\) for
\(\ell\ge2\), the bound~\eqref{eq:degree} gives
\[
  B(C)\le2e(\mathcal H')\le2L^{r-1}e(C)\,\Delta_2(\mathcal H')
  =O\bigl(L^{r-1}\,2^{-(s-1)L}n^2\bigr)=o(n^2).\qedhere
\]
\end{proof}

The deferred copies enter through the relation system. For
\(C\in\mathcal C\) and each row \(a\), let \(C_a\) be the set of edges of
\(C\) at \(a\) (with \(a\) as the row), and let \(\Gamma_a\) be the graph
on \(C_a\) in which \(ab\) and \(ab'\) are adjacent if \(b\) and \(b'\)
are partners, that is, if \(bb'\in D\). Thus \(\Gamma_a\) is a copy of
\(D[N_C(a)]\). Every edge belongs to at most two sets \(C_a\), so
\(\sum_a|C_a|\le 2|C|\). We partition each \(\Gamma_a\) as in
Section~\ref{sec:deferred}, using a fixed ordering, and identify its
parts with sets of columns. Let \(\mathcal P_a(C)\) be the family of
small parts at row \(a\).
Then \(\mathcal T(C)\) consists of the subgraphs of \(C\) that contain no
edge of a large part and, for every row \(a\) and every
\(P\in\mathcal P_a(C)\), contain at most one edge \(ab\) with \(b\in P\).
By~(P2), two columns in distinct small parts at the same row are not
partners.

\begin{lemma}\label{lem:kst-short}
For every \(C\in\mathcal C\), every admissible \(I\subseteq C\)
contains at most \(M=(t-1)n\,2^{2L^2}=n^{1+o(1)}\) short pairs.
\end{lemma}

\begin{proof}
A short pair contained in \(I\) has the form \(\{ab,ab'\}\), where the
columns \(b,b'\) are at distance less than \(L\) in \(D\). Since
\(\Delta(D)<2^{2L}\), at most \(n\,2^{2L^2}\) unordered pairs of columns
are at distance less than \(L\) in \(D\), and by (A2) each of them has at
most \(t-1\) common rows in \(I\).
\end{proof}

Condition (A2) is used only here and in the proof of
Lemma~\ref{lem:kst-columns}.

\begin{lemma}\label{lem:kst-reduction}
We have
\[
  |\mathcal A(G)|\le2^{o(n^{1+1/s})}\max_{C\in\mathcal C}
  |\mathcal T(C)|.
\]
\end{lemma}

\begin{proof}
Since \(L\to\infty\), we have \(L\ge8r^2\) for large \(n\), and by
Lemma~\ref{lem:kst-short}, Theorem~\ref{thm:deferred} applies with
\(M=n^{1+o(1)}\). Since \(e(G_0)=O(n^{1+1/s})\), the exponential factor
in~\eqref{eq:deferred-count} is \(2^{o(n^{1+1/s})}\). Combine this with
Lemma~\ref{lem:kst-columns}(c).
\end{proof}

\subsection{A weighted K\H{o}v\'ari--S\'os--Tur\'an bound}\label{sec:kst-weights}

By Lemma~\ref{lem:kst-reduction}, it remains to bound
\(|\mathcal T(C)|\) for a single container \(C\in\mathcal C\). A member of
\(\mathcal T(C)\) uses at most one column of each small part \(P\) at
each row, so a small part behaves like a single column with \(1+|P|\)
possible choices. We give it the weight \(|P|^{1-1/s}\), which is
at least \(\log(1+|P|)\) when \(s\ge3\), and let
\begin{equation}\label{eq:kst-weight}
  \Phi(C)=\sum_{\text{rows }a}\ \sum_{P\in\mathcal P_a(C)}|P|^{1-1/s}.
\end{equation}
If all small parts were single columns, \(\Phi(C)\) would be the number
of edges in them, and Lemma~\ref{lem:kst} would bound it when \(C\) is
admissible. The next lemma is a weighted form of this bound, in which
the packings of Lemma~\ref{lem:kst-packings} account for the \(s\)-sets
of rows with too many common columns.

\begin{lemma}\label{lem:host-records}
Every \(C\in\mathcal C\) satisfies
\[
  \Phi(C)\le\mu^{1/s}mn^{1-1/s}+o(n^{1+1/s}).
\]
\end{lemma}

\begin{proof}
We repeat the proof of Lemma~\ref{lem:kst} with weights. If a column
\(b\) lies in a small part \(P\) of the row \(a\), let
\(w_a(b)=|P|^{-1/s}\), and otherwise let \(w_a(b)=0\). A small part \(P\)
has total weight \(|P|\cdot|P|^{-1/s}\) in its row, so
\(\Phi(C)=\sum_bw(b)\), where \(w(b)=\sum_aw_a(b)\). For an \(s\)-set \(U\)
of rows, let \(\rho(U)=\sum_b\prod_{a\in U}w_a(b)\); with unit weights,
this would be the number of common columns of \(U\).

\begin{claim}\label{clm:kst-rho}
For every \(s\)-set \(U\) of rows,
\(\rho(U)\le\mu+(\mu+1)|\mathcal B_U(C)|\).
\end{claim}

\begin{proof}
Group the columns \(b\) with \(\prod_{a\in U}w_a(b)>0\) according to the
tuple \((P_a)_{a\in U}\) of small parts containing them. A group
contributes \(|\bigcap_aP_a|/\prod_a|P_a|^{1/s}\le1\) to \(\rho(U)\). Two
columns in distinct groups lie in distinct small parts of some row
\(a\in U\), so they are not partners. Hence, after the
\((\mu+1)|\mathcal B_U(C)|\) columns of the members of
\(\mathcal B_U(C)\) are deleted, at most \(\mu\) groups remain nonempty,
since otherwise one column from each of \(\mu+1\) groups would extend
\(\mathcal B_U(C)\). Each deleted column contributes at most one to
\(\rho(U)\).
\end{proof}

As in Lemma~\ref{lem:kst}, we separate the columns of small weight. Call
a column \(b\) \emph{heavy} if \(w(b)\ge s-1\); the other columns
contribute less than \((s-1)n\) to \(\Phi(C)\). Let \(b\) be heavy, and
choose distinct rows \(a_1,\dots,a_s\) one at a time. Since
\(w_a(b)\le1\) for every row \(a\), when \(a_j\) is chosen the remaining
rows have total weight at least \(w(b)-j+1\ge w(b)-s+1\ge0\) at \(b\).
Hence the sum of \(\prod_jw_{a_j}(b)\) over all ordered choices is at
least \((w(b)-s+1)^s\). Summing over heavy columns and using
Claim~\ref{clm:kst-rho}, we obtain
\[
  \sum_{b\text{ heavy}}\bigl(w(b)-s+1\bigr)^s\le s!\sum_U\rho(U)
  \le\mu m^s+s!(\mu+1)B(C).
\]
Since there are at most \(n\) columns, H\"older's inequality gives
\[
  \Phi(C)\le(s-1)n+\sum_{b\text{ heavy}}\bigl(w(b)-s+1\bigr)
  \le(s-1)n+n^{1-1/s}\bigl(\mu m^s+s!(\mu+1)B(C)\bigr)^{1/s}.
\]
By Lemma~\ref{lem:kst-packings}, \(B(C)=o(n^2)\), and since
\(m\le\lambda n^{2/s}\), the right-hand side is
\(\mu^{1/s}mn^{1-1/s}+o(n^{1+1/s})\).
\end{proof}

For \(s\ge3\), Lemma~\ref{lem:host-records} suffices to complete the
proof.

\begin{proof}[Proof of Proposition~\ref{prop:local} for \(s\ge3\)]
Fix \(C\in\mathcal C\). Since the host is bipartite, the sets \(C_a\) are
disjoint, and~\eqref{eq:product} gives
\[
  |\mathcal T(C)|
  =\prod_a\prod_{P\in\mathcal P_a(C)}(1+|P|).
\]
For a positive integer \(x\), we have
\(\log(1+x)\le x^{2/3}\le x^{1-1/s}\). The first inequality holds for
\(x=1,2\) by direct computation, and hence for all \(x\ge2\), because
\(x^{-2/3}\log(1+x)\) is decreasing on \([2,\infty)\). By
Lemma~\ref{lem:host-records},
\[
  \log|\mathcal T(C)|\le\Phi(C)
  \le\mu^{1/s}mn^{1-1/s}+o(n^{1+1/s}),
\]
and Lemma~\ref{lem:kst-reduction} completes the proof.
\end{proof}

For \(s=2\), bounding~\eqref{eq:product} part by part would lose a factor
of about two, because an edge \(ab\) of a member of \(\mathcal T(C)\)
occupies a small part at the row \(a\) and another at the row \(b\). We
therefore regard the members of \(\mathcal T(C)\) as matchings in an
auxiliary graph, and the factor \(\frac12\) comes from the following
bound of~\cite{yang}. For a graph \(H\), let \(Z(H)\) be the number of
matchings of \(H\), including the empty matching.

\begin{proposition}[\cite{yang}]\label{prop:matchings}
Every graph \(H\) satisfies
\[
  \log Z(H)\le\tfrac12\sum_{v\in V(H)}\sqrt{d_H(v)}.
\]
\end{proposition}

\begin{proof}[Proof of Proposition~\ref{prop:local} for \(s=2\)]
Here \(G\) is a graph on \([n]\) and \(m=n\). For \(C\in\mathcal C\),
form a graph \(\widehat C\) with one vertex \((a,P)\) for each row \(a\)
and each \(P\in\mathcal P_a(C)\). If an edge \(ab\) of \(C\) lies in no
large part, then \(b\) lies in a small part \(P\) of the row \(a\) and
\(a\) lies in a small part \(P'\) of the row \(b\); we join \((a,P)\) to
\((b,P')\) in \(\widehat C\). This correspondence identifies the sets in
\(\mathcal T(C)\) with the matchings of \(\widehat C\). Since
\(d_{\widehat C}(a,P)\le|P|\), Proposition~\ref{prop:matchings} and
Lemma~\ref{lem:host-records} give
\[
  \log|\mathcal T(C)|
  =\log Z(\widehat C)
  \le\tfrac12\sum_a\sum_{P\in\mathcal P_a(C)}\sqrt{|P|}
  =\tfrac12\Phi(C)
  \le\tfrac12\sqrt\mu\,n^{3/2}+o(n^{3/2}).
\]
Lemma~\ref{lem:kst-reduction} completes the proof.
\end{proof}

\subsection{Projections}\label{sec:kst-projections}

In this subsection \(s\ge3\). Fix \(\lambda>0\) and a graph \(G\) on
\([n]\) with \(\Delta(G)\le\lambda n^{1-1/s}\). The rows of each
bipartite host are taken from the common neighbourhood of an
\((s-2)\)-set, and Proposition~\ref{prop:local} bounds the number of
admissible subgraphs of each host. Before applying Shearer's inequality,
we remove the vertices of small degree from the graph being counted.
The hosts can then be chosen so that they cover almost every remaining
edge twice, and this halves the sum of the local bounds. For \(W\subseteq[n]\) and \(d>0\), let
\(\mathcal I_d(G[W])\) be
the family of those \(I\in\mathcal I(G[W])\) in which every vertex of
\(W\) has degree at least \(d\).

\begin{lemma}\label{lem:projections}
Let \(\varepsilon>0\), and let \(d_0=\lceil\varepsilon n^{1-1/s}\rceil\).
For every \(W\subseteq[n]\),
\[
  |\mathcal I_{d_0}(G[W])|\le2^{c_{s,t}|W|n^{1-1/s}+o(n^{2-1/s})},
\]
where the error term depends on \(\lambda\) and \(\varepsilon\) but not
on \(G\) and \(W\).
\end{lemma}

We first deduce Proposition~\ref{prop:bounded}.

\begin{proof}[Proof of Proposition~\ref{prop:bounded} for \(s\ge3\),
assuming Lemma~\ref{lem:projections}]
We remove the vertices of small degree and apply
Lemma~\ref{lem:projections} to the rest. Choose
\(0<\varepsilon<\lambda\) with
\(\varepsilon\log(e\lambda/\varepsilon)\le c_{s,t}\), and let
\(d_0=\lceil\varepsilon n^{1-1/s}\rceil\). For \(I\in\mathcal I(G)\),
repeatedly delete a vertex of degree less than \(d_0\) in what remains
of \(I\), chosen according to a fixed ordering of \([n]\), until every
remaining vertex has degree at least \(d_0\); let \(W\) be the set of
remaining vertices. Let \(\mathcal D_W\) be the family of those
\(I\in\mathcal I(G)\) that yield \(W\). Then
\(I\cap E(G[W])\in\mathcal I_{d_0}(G[W])\) for every \(I\in\mathcal
D_W\).

The edges of \(I\) outside \(G[W]\) are the edges deleted in the
process. There are at most \(\varepsilon(n-|W|)n^{1-1/s}\) of them, and
they lie among the at most \(\lambda(n-|W|)n^{1-1/s}\) edges of \(G\)
meeting \([n]\setminus W\). By~\eqref{eq:binomial} and the choice of
\(\varepsilon\), the set \(I\setminus E(G[W])\) takes at most
\(2^{c_{s,t}(n-|W|)n^{1-1/s}}\) values. By Lemma~\ref{lem:projections},
\[
  |\mathcal I(G)|=\sum_{W\subseteq[n]}|\mathcal D_W|
  \le\sum_{W\subseteq[n]}2^{c_{s,t}(n-|W|)n^{1-1/s}}
     \bigl|\mathcal I_{d_0}(G[W])\bigr|
  \le2^n\,2^{c_{s,t}n^{2-1/s}+o(n^{2-1/s})},
\]
and \(n=o(n^{2-1/s})\).
\end{proof}

For a family \(\mathcal K\) of edge sets and an edge set \(E\), we write
\(\mathcal K|_E=\{K\cap E:K\in\mathcal K\}\). We use the projection
inequality of Chung, Graham, Frankl and Shearer~\cite{shearer} in the
following form: if every element of a finite set \(E\) lies in at least
two of the sets \(E_1,\dots,E_q\subseteq E\), then every family
\(\mathcal K\) of subsets of \(E\) satisfies
\begin{equation}\label{eq:shearer}
  2\log|\mathcal K|\le\sum_{i=1}^q\log\bigl|\mathcal K|_{E_i}\bigr|.
\end{equation}

\begin{proof}[Proof of Lemma~\ref{lem:projections}]
We may assume that \(W\neq\varnothing\). To each
\(I\in\mathcal I_{d_0}(G[W])\) we assign a sequence of \((s-2)\)-sets,
together with disjoint sets of rows taken from their common
neighbourhoods in \(I\), such that the corresponding bipartite hosts
cover almost every edge of \(I\) twice. For a fixed sequence,
Proposition~\ref{prop:local} bounds the number of possible restrictions
of \(I\) to each host, and~\eqref{eq:shearer} combines these bounds. We regard each \(I\) as a
graph on \([n]\), with the vertices outside \(W\) isolated, and begin by
bounding its common neighbourhoods.

\begin{claim}\label{clm:common}
Let \(I\in\mathcal I_{d_0}(G[W])\) and \(X\subseteq[n]\) with
\(1\le|X|=j<s\). Then \(|N_I(X)|=O(n^{1-j/s})\), where the implied
constant depends on \(\varepsilon\).
\end{claim}

\begin{proof}
We may assume that \(X\subseteq W\). We count pairs \((z,Y)\) with
\(z\in N_I(X)\) and \(Y\) an \((s-j)\)-subset of
\(N_I(z)\setminus X\). Each \(z\) lies in at least
\(\binom{d_0-j}{s-j}\) pairs. For each \(Y\) there are at most \(t-1\)
choices of \(z\), since the \(s\)-set \(X\cup Y\) has at most \(t-1\)
common neighbours in \(I\). Hence
\(|N_I(X)|\binom{d_0-j}{s-j}\le(t-1)\binom{n-j}{s-j}\), and the claim
follows from \(d_0\ge\varepsilon n^{1-1/s}\).
\end{proof}

Let \(\eta>0\). Choose \(0<\alpha<1\) with
\(\lambda\alpha\log(e/\alpha)\le\eta\), then \(\beta\ge1\) with
\(2^{-\beta}<\alpha/2\), and let
\[
  q=\left\lceil\beta\binom n{s-2}\Bigm/\binom{d_0}{s-2}\right\rceil
    =O\bigl(n^{1-2/s}\bigr).
\]
A \emph{projection record} is a sequence \(\pi=(X_i,Y_i)_{i=1}^q\) in
which each \(X_i\) is an \((s-2)\)-subset of \([n]\) and the \(Y_i\) are
disjoint subsets of \(W\). Given \(\pi\), let \(E_i\) be the set of edges
of \(G[W]\) with an endpoint in \(Y_i\) and no endpoint in \(X_i\), and
let \(E^*\) be the set of edges lying in two of the sets \(E_i\). Since
the \(Y_i\) are disjoint, every edge lies in at most two of the sets
\(E_i\). There are at most
\(\binom n{s-2}^q(q+1)^n=2^{o(n^{2-1/s})}\) projection records.

\begin{claim}\label{clm:coverage}
Every \(I\in\mathcal I_{d_0}(G[W])\) has a projection record with
\(Y_i\subseteq N_I(X_i)\) for every \(i\) and
\(|I\setminus E^*|\le\alpha\lambda n^{2-1/s}\).
\end{claim}

\begin{proof}
Choose \(X_1,\dots,X_q\) independently and uniformly at random. Assign
each \(v\in W\) to \(Y_i\) for the first index \(i\) with
\(X_i\subseteq N_I(v)\), and leave \(v\) unassigned if there is no such
index. Then \(Y_i\subseteq N_I(X_i)\). An edge \(uv\) of \(I\) lies in
\(E^*\) if \(u\) and \(v\) are assigned to distinct indices and each of
them lies outside the set \(X_i\) to which the other is assigned. We
bound the probability that one of these conditions fails.

Each endpoint is unassigned with probability at most
\[
  \left(1-\binom{d_0}{s-2}\Bigm/\binom n{s-2}\right)^q\le 2^{-\beta},
\]
since its neighbourhood in \(I\) has at least \(d_0\) elements.
Conditional on \(v\) being assigned, the set \(X_i\) to which it is
assigned is uniformly distributed among the \((s-2)\)-subsets of
\(N_I(v)\), so it contains \(u\) with probability at most
\((s-2)/d_0=o(1)\). Finally, if \(u\) and \(v\) are assigned to the same
index, then the first set \(X_i\) contained in \(N_I(u)\) or in
\(N_I(v)\) is contained in \(N_I(\{u,v\})\). Conditional on such a set
existing, it is uniformly distributed among the \((s-2)\)-sets contained
in \(N_I(u)\) or in \(N_I(v)\), so this happens with probability at most
\(\binom{|N_I(\{u,v\})|}{s-2}/\binom{d_0}{s-2}=o(1)\), by
Claim~\ref{clm:common} with \(j=2\).

The expected number of edges of \(I\) outside \(E^*\) is therefore at
most \((2^{1-\beta}+o(1))e(I)\). Since \(e(I)\le\frac\lambda2n^{2-1/s}\)
and \(2^{-\beta}<\alpha/2\), some outcome satisfies the claim for all
sufficiently large \(n\).
\end{proof}

Choose one projection record as in Claim~\ref{clm:coverage} for each
\(I\in\mathcal I_{d_0}(G[W])\). Let \(\mathcal R\) be the set of chosen
records, and let \(\mathcal K_\pi\) be the family of those \(I\) for which
\(\pi\) was chosen.

\begin{claim}\label{clm:projections}
Let \(\pi\in\mathcal R\) and \(\mathcal K^*=\mathcal K_\pi|_{E^*}\). For
every \(i\),
\[
  \bigl|\mathcal K^*|_{E_i}\bigr|\le
  2^{\mu^{1/s}|Y_i|n^{1-1/s}+o(n^{1+1/s})},
\]
where the error term does not depend on \(\pi\) or \(i\).
\end{claim}

\begin{proof}
Let \(G_i\) be the host with rows \(Y_i\), columns \(W\setminus X_i\),
and an edge between \(a\) and \(b\) whenever \(ab\in E(G)\). Every edge
in \(E_i\cap E^*\) also lies in some \(E_j\) with \(j\ne i\). Since the
row sets \(Y_i\) and \(Y_j\) are disjoint, this edge has exactly one
endpoint in \(Y_i\). Thus there is an injection
\[
  \phi_i\colon E_i\cap E^*\longrightarrow E(G_i),\qquad
  \{a,b\}\longmapsto(a,b)
  \quad(a\in Y_i,\ b\in W\setminus X_i).
\]
For \(I\in\mathcal K_\pi\), the image \(\phi_i(I\cap E_i\cap E^*)\) is a
spanning subgraph of \(G_i\). Its rows \(Y_i\) lie in \(N_I(X_i)\) and
its columns \(W\setminus X_i\) lie outside \(X_i\). It is therefore a
subgraph of the admissible bipartite graph described in
Section~\ref{sec:kst-hosts} for \(I\) and \(X_i\), so it belongs to
\(\mathcal A(G_i)\). As \(\phi_i\) is injective,
\(|\mathcal K^*|_{E_i}|\le|\mathcal A(G_i)|\).

By Claim~\ref{clm:common} with \(j=s-2\), there is a constant
\(\lambda'\ge\lambda\), depending only on \(\lambda\) and
\(\varepsilon\), such that \(G_i\) has at most \(\lambda'n^{2/s}\) rows,
each of degree at most \(\lambda'n^{1-1/s}\). Proposition~\ref{prop:local}
with \(\lambda'\) in place of \(\lambda\) therefore gives
\(\log|\mathcal A(G_i)|\le\mu^{1/s}|Y_i|n^{1-1/s}+o(n^{1+1/s})\), with
the error uniform in \(\pi\) and \(i\).
\end{proof}

Let \(\pi\in\mathcal R\) and
\(\mathcal K^*=\mathcal K_\pi|_{E^*}\). On the ground set \(E^*\), each
edge lies in exactly two of the sets \(E_i\cap E^*\), and
\(\mathcal K^*|_{E_i\cap E^*}=\mathcal K^*|_{E_i}\).
By~\eqref{eq:shearer} and Claim~\ref{clm:projections},
\[
  \begin{aligned}
    2\log|\mathcal K^*|
    &\le\sum_{i=1}^q\log\bigl|\mathcal K^*|_{E_i}\bigr|\\
    &\le\mu^{1/s}n^{1-1/s}\sum_{i=1}^q|Y_i|+o(qn^{1+1/s})\\
    &\le\mu^{1/s}|W|n^{1-1/s}+o(n^{2-1/s}),
  \end{aligned}
\]
where the last step uses \(\sum_i|Y_i|\le|W|\) and \(q=O(n^{1-2/s})\).
On the other hand, \(I\setminus E^*\) is a subset of \(E(G[W])\) with at
most \(\alpha\lambda n^{2-1/s}\) elements, and
\(e(G[W])\le\lambda n^{2-1/s}\), so by~\eqref{eq:binomial} it takes at
most \((e/\alpha)^{\alpha\lambda n^{2-1/s}}\le2^{\eta n^{2-1/s}}\)
values. Since \(I\) is determined by \(I\cap E^*\) and \(I\setminus E^*\),
and since \(\mu^{1/s}/2=c_{s,t}\), we conclude that
\[
  |\mathcal I_{d_0}(G[W])|=\sum_{\pi\in\mathcal R}|\mathcal K_\pi|
  \le|\mathcal R|\,2^{\eta n^{2-1/s}}\max_{\pi\in\mathcal R}
     \bigl|\mathcal K_\pi|_{E^*}\bigr|
  \le2^{c_{s,t}|W|n^{1-1/s}+(\eta+o(1))n^{2-1/s}}.
\]
Since \(\eta>0\) was arbitrary, the lemma follows.
\end{proof}

\section{Sidon sets}\label{sec:sidon}

In this section we prove Theorem~\ref{thm:sidon}. For \(U\subseteq[n]\),
let \(\mathcal A(U)\) be the family of Sidon subsets of \(U\), so that
\(F(n)=|\mathcal A([n])|\). Throughout this section \(n\to\infty\).

For completeness, we recall the lower-bound construction of Saxton and
Thomason~\cite{stonline}. Let \(q\) be the largest prime with
\(4q(q-1)\le n\), and take Ruzsa's Sidon set of size \(q-1\) modulo
\(q(q-1)\)~\cite{ruzsa}. Choosing either none or one of the four
representatives in \([4q(q-1)]\) of each element of this set gives
\(5^{q-1}\) Sidon sets, and since \(q=(\frac12+o(1))\sqrt n\), we have
\(F(n)\ge5^{q-1}=2^{(\frac12\log 5+o(1))\sqrt n}\).

For the upper bound, we bound \(|\mathcal A(U)|\) while the set \(U\)
is successively replaced by smaller sets. The argument follows the same
steps as in Section~\ref{sec:kst}, and corresponding objects are denoted
by the same letters: \(\Hact\) for the active configurations,
\(\mathcal C\) for the containers, \(\Gamma\) for the relation graph,
\(P\) for a small part, \(w\) for the weights and \(\Phi(C)\) for the
total weight. In Section~\ref{sec:sidon-coarse} an iteration of the
container lemma reduces the problem to a single coarse container \(U\)
with \(|U|=O(\sqrt n)\) spanning \(O(n)\) additive quadruples
(Lemma~\ref{lem:sidon-first-pass}). The main estimate of the section,
Proposition~\ref{prop:terminal}, bounds \(|\mathcal A(U)|\) for such a
set. The standard iteration cannot be continued beyond this point,
since some sums and differences may still have about \(\sqrt n\)
representations in \(U\). In Section~\ref{sec:sidon-popular} we defer
the quadruples with such a sum or difference, and in
Section~\ref{sec:sidon-containers} we apply Theorem~\ref{thm:deferred} to
the remaining ones. This gives
\[
  |\mathcal A(U)|\le2^{o(\sqrt n)}\max_{C}\prod_{P}\bigl(1+|P|\bigr)
\]
(Lemma~\ref{lem:sidon-reduction}), where \(C\) runs over the containers
and \(P\) over the small parts of \(C\). Finally, in
Section~\ref{sec:sidon-et} a weighted form of the Erd\H{o}s--Tur\'an
argument shows that \(\sum_P\sqrt{|P|}\le(1+o(1))\sqrt n\), and the
inequality \(\log(1+a)\le\frac12\log5\cdot\sqrt a\) completes the proof.

\subsection{Reduction to a coarse container}\label{sec:sidon-coarse}

An \emph{additive quadruple} in \(U\) is a tuple \((u,v,w,x)\in U^4\)
with \(u<v\le w<x\) and \(u+x=v+w\). Its sum is \(u+x\), and its
differences are \(v-u=x-w\) and \(w-u=x-v\). The quadruple is
\emph{degenerate} if \(v=w\), in which case \(u,v,x\) form an arithmetic
progression. A set \(A\subseteq U\) is a Sidon set if and only if it
contains no additive quadruple: if \(a+b=c+d\) with \(a\le b\),
\(c\le d\) and \(\{a,b\}\neq\{c,d\}\), then the two pairs are disjoint,
and the two extreme and the two middle elements of the multiset
\(\{a,b,c,d\}\) form the two pairs. A Sidon set also has at most one
pair of elements at each positive distance, since \(x-y=x'-y'\) for two
different pairs with \(x>y\) and \(x'>y'\) gives the repeated sum
\(x+y'=x'+y\).

Let \(\mathcal Q(U)\) be the \(4\)-uniform hypergraph on \(U\) whose edges
are the sets of entries of the nondegenerate additive quadruples, and let
\(Q(U)=e(\mathcal Q(U))\). For \(h\ge1\), let \(r_U^-(h)\) be the number
of pairs of elements of \(U\) at distance \(h\), and for \(m\in\mathbb Z\),
let \(r_U^+(m)\) be the number of pairs of distinct elements of \(U\)
with sum \(m\). Two distinct pairs with the same sum form a nondegenerate
quadruple. Two distinct pairs at the same distance form a nondegenerate
quadruple if they are disjoint, and a degenerate one if they meet; a
nondegenerate quadruple arises in this way at both of its differences.
The degenerate quadruples with middle entry \(a\) are counted by
\(r_U^+(2a)\). Therefore
\begin{align}\label{eq:representations}
  Q(U)&=\sum_m\binom{r_U^+(m)}2,\\
  \sum_{h\ge1}\binom{r_U^-(h)}2
    &=2Q(U)+\sum_{a\in U}r_U^+(2a),\notag
\end{align}
and, by the Cauchy--Schwarz inequality,
\begin{equation}\label{eq:progressions}
  \sum_{a\in U}r_U^+(2a)
  \le\sqrt{|U|\sum_m r_U^+(m)^2}
  =\sqrt{|U|\left(2Q(U)+\binom{|U|}2\right)}.
\end{equation}
The last equality uses the first identity
in~\eqref{eq:representations} and \(\sum_m r_U^+(m)=\binom{|U|}2\).

Two entries of a nondegenerate quadruple, together with the other two,
form either two pairs with the same sum or two pairs at the same
distance, and at most \(|U|\) pairs have a given sum or a given distance.
Moreover, three entries determine the fourth in at most three ways,
according to which of them is paired with it. Hence
\begin{equation}\label{eq:sidon-codegrees}
  \Delta_2(\mathcal Q(U))\le2|U|,\qquad
  \Delta_3(\mathcal Q(U))\le3,\qquad
  \Delta_4(\mathcal Q(U))\le1.
\end{equation}
The first container step is a standard iteration of
Theorem~\ref{thm:standard}.

\begin{lemma}\label{lem:sidon-first-pass}
For every \(\eta>0\) there is \(\lambda>0\) such that, for all
sufficiently large \(n\), there is a family \(\mathcal U\) of at most
\(2^{\eta\sqrt n/2}\) subsets of \([n]\) with
\[
  |U|\le\lambda\sqrt n,\quad Q(U)\le\lambda n\quad(U\in\mathcal U),
  \qquad \mathcal A([n])=\bigcup_{U\in\mathcal U}\mathcal A(U).
\]
\end{lemma}

\begin{proof}
We first show that a set spanning few quadruples is small. Suppose that
\(Q(U)\le xn\) with \(x\ge1\). There are fewer than \(2n\) possible
sums, and \(\sum_m r_U^+(m)=\binom{|U|}2\), so the Cauchy--Schwarz
inequality and~\eqref{eq:representations} give
\[
  \frac{1}{4n}\binom{|U|}{2}^{\!2}-\frac12\binom{|U|}{2}
  \le Q(U)\le xn.
\]
It follows that \(\binom{|U|}{2}\le n(1+\sqrt{1+4x})\), and so
\begin{equation}\label{eq:size-from-quadruples}
  |U|=O(\sqrt n\,x^{1/4}).
\end{equation}

We build families \(\mathcal U_0,\mathcal U_1,\dots\) in stages,
maintaining \(Q(U)\le x_jn\) for every \(U\in\mathcal U_j\), where
\(x_j=2^{-j}n^2\). Let \(c_0,\beta\ge1\) be constants to be chosen.
Start with \(\mathcal U_0=\{[n]\}\); since the three smallest entries
determine a nondegenerate quadruple, \(Q([n])\le x_0n\). We stop at the
first \(j\) with \(x_j<\beta\). If \(x_j\ge\beta\), we put each
\(U\in\mathcal U_j\) with \(Q(U)\le x_{j+1}n\) into
\(\mathcal U_{j+1}\); for every other \(U\in\mathcal U_j\), we apply
Theorem~\ref{thm:standard} to \(\mathcal Q(U)\) with parameter
\(p_j=c_0x_j^{-1/3}\) and put the resulting containers into
\(\mathcal U_{j+1}\). By~\eqref{eq:sidon-codegrees}
and~\eqref{eq:size-from-quadruples}, if \(Q(U)>x_jn/2\), then
\[
  \frac{|U|}{Q(U)}\max_{2\le i\le4}p_j^{1-i}\Delta_i(\mathcal Q(U))
  =O\Bigl(\frac{x_j^{-1/6}}{c_0}+\frac{x_j^{-1/12}}{c_0^2\sqrt n}
    +\frac{x_j^{1/4}}{c_0^3\sqrt n}\Bigr).
\]
As \(1\le x_j\le n^2\), we may choose \(c_0\) so that the right-hand side
is at most \(1/2\), and then \(\beta\) so large that \(p_j\le1/128\).
By~\eqref{eq:degree}, every resulting container \(C\subseteq U\)
satisfies \(Q(C)\le Q(U)/2\le x_{j+1}n\), and every Sidon subset of \(U\)
lies in one of them. Hence
\(\bigcup_{U\in\mathcal U_j}\mathcal A(U)=\mathcal A([n])\) at every
stage.

By~\eqref{eq:fingerprints} and~\eqref{eq:size-from-quadruples}, the
number of containers produced from one \(U\) at stage \(j\) is at most
\(2^{O(p_j|U|\log(e/p_j))}=2^{O(\sqrt n\,x_j^{-1/12}\log x_j)}\). The sum of
these exponents over the stages with \(x_j\ge\beta\) is
\(O(\sqrt n\,\beta^{-1/12}\log\beta)\), so if \(\beta\) is large enough,
the final family has at most \(2^{\eta\sqrt n/2}\) members. Each
member \(U\) of the final family satisfies \(Q(U)<\beta n\),
and~\eqref{eq:size-from-quadruples} bounds \(|U|\).
\end{proof}

The main estimate of the section counts the Sidon sets in a single
coarse container.

\begin{proposition}\label{prop:terminal}
Let \(\lambda>0\) be fixed. Every \(U\subseteq[n]\) with
\(|U|\le\lambda\sqrt n\) and \(Q(U)\le\lambda n\) satisfies
\[
  \log|\mathcal A(U)|\le\bigl(\tfrac12\log 5+o(1)\bigr)\sqrt n,
\]
where the \(o(1)\) term depends on \(\lambda\) but not on \(U\).
\end{proposition}

\begin{proof}[Proof of Theorem~\ref{thm:sidon} assuming
Proposition~\ref{prop:terminal}]
The lower bound was shown above. For the upper bound, let
\(\eta>0\), and take the family \(\mathcal U\) and the constant
\(\lambda\) from Lemma~\ref{lem:sidon-first-pass}.
Proposition~\ref{prop:terminal} gives
\[
  F(n)\le|\mathcal U|\max_{U\in\mathcal U}|\mathcal A(U)|
  \le2^{(\frac12\log 5+\eta/2+o(1))\sqrt n}.
\]
Letting \(n\to\infty\) and then \(\eta\to0\) gives the upper bound.
\end{proof}

\subsection{Popular sums and differences}\label{sec:sidon-popular}

For the rest of this section, fix \(\lambda>0\) and \(U\subseteq[n]\) with
\(|U|\le\lambda\sqrt n\) and \(Q(U)\le\lambda n\). Implied constants and
\(o(\cdot)\) terms may depend on \(\lambda\) but not on \(U\). As in
Section~\ref{sec:kst}, we set \(L=(\log n)^{1/3}\), treated as an
integer.

The container iteration cannot be continued because
\(\Delta_2(\mathcal Q(U))\) may be of order \(\sqrt n\): if a sum \(m\) has
\(r_U^+(m)\) representations, then two elements with sum \(m\) lie
together in \(r_U^+(m)-1\) quadruples. On the other hand, a Sidon set
contains at most one pair with a given sum and at most one pair at a
given distance. The quadruples through a sum or a difference with many
representations can therefore be used inside the containers instead. We call a sum
\(m\in\mathbb Z\) or a difference \(h\ge1\) \emph{popular} if it has at
least \(\sqrt n/\log n\) representations in \(U\), and let
\[
  D^+=\Bigl\{m\in\mathbb Z:r_U^+(m)\ge\frac{\sqrt n}{\log n}\Bigr\},\qquad
  D^-=\Bigl\{h\ge1:r_U^-(h)\ge\frac{\sqrt n}{\log n}\Bigr\}.
\]
Let \(\mathcal H=\mathcal Q(U)\). We defer an edge of \(\mathcal H\) if
its sum belongs to \(D^+\) or one of its two differences belongs to
\(D^-\), and let \(\Hact\) consist of the other edges, which we call
\emph{active}. For \(C\subseteq U\), write \(\Qact(C)=e(\Hact[C])\).
Every Sidon subset of \(U\) is independent in \(\mathcal H\), and hence in
\(\Hact\).

\begin{lemma}\label{lem:sidon-active}
The following hold.
\begin{statements}
\item \(|D^+|+|D^-|=O\bigl((\log n)^2\bigr)\).
\item \(\Delta_2(\Hact)\le2\sqrt n/\log n\).
\end{statements}
\end{lemma}

\begin{proof}
By~\eqref{eq:progressions} and the bounds on \(|U|\) and \(Q(U)\), we
have \(\sum_ar_U^+(2a)=O(n^{3/4})\), so both binomial sums
in~\eqref{eq:representations} are \(O(n)\). Each popular value
contributes \(\Omega(n/(\log n)^2)\) to one of these sums, which
proves~(a).

For (b), let \(x<y\) be two entries of an active quadruple. The other two
entries form a pair with sum \(x+y\) or a pair at distance \(y-x\), and
since the quadruple is active, this sum or distance is not popular. So
there are fewer than \(2\sqrt n/\log n\) possibilities for the other two
entries.
\end{proof}

\subsection{Containers and small parts}\label{sec:sidon-containers}

We apply Theorem~\ref{thm:deferred} to \(\Hact\) on \(U\), with
\(\mathcal F=\mathcal A(U)\) and \(r=4\), and let \(\mathcal C\) be the
resulting family of containers. Since the active quadruples have small
codegrees by Lemma~\ref{lem:sidon-active}(b), each container spans few
of them.

\begin{lemma}\label{lem:sidon-containers}
For all sufficiently large \(n\), every \(C\in\mathcal C\) satisfies
\begin{equation}\label{eq:active-few}
  \Qact(C)\le\frac{2\lambda n}{(\log n)^{2/3}}.
\end{equation}
\end{lemma}

\begin{proof}
Apply~\eqref{eq:degree} to \(\Hact[C]\), with \(p^{1-\ell}=L^{\ell-1}\).
By Lemma~\ref{lem:sidon-active}(b) and~\eqref{eq:sidon-codegrees}, whose
bounds hold for every subhypergraph of \(\mathcal Q(U)\),
\[
  \Qact(C)
  \le |C|\max\Bigl\{\frac{2L\sqrt n}{\log n},\,3L^2,\,L^3\Bigr\}
  \le\frac{2\lambda n}{(\log n)^{2/3}}
\]
for sufficiently large \(n\), since \(|C|\le\lambda\sqrt n\) and
\(L^3=\log n\).
\end{proof}

The deferred quadruples enter through a single relation graph. For
\(C\in\mathcal C\), let \(\Gamma\) be the graph on \(C\) in which two
distinct elements \(x,y\) are adjacent if \(x+y\in D^+\) or
\(|x-y|\in D^-\). We partition \(\Gamma\) as in Section~\ref{sec:deferred},
using the natural order. Let \(\mathcal P(C)\) be the family of small
parts. By~\eqref{eq:product},
\begin{equation}\label{eq:sidon-product}
  |\mathcal T(C)|=\prod_{P\in\mathcal P(C)}(1+|P|).
\end{equation}
The short-pair hypothesis of Theorem~\ref{thm:deferred} holds because a
path in \(\Gamma\) is a composition of reflections and translations,
and a Sidon set realises each sum and each difference at most once.

\begin{lemma}\label{lem:sidon-short}
For every \(C\in\mathcal C\), every Sidon set \(A\subseteq C\) contains
at most \(M=(1+|D^+|+2|D^-|)^{L-1}=n^{o(1)}\) short pairs.
\end{lemma}

\begin{proof}
For each short pair in \(A\), fix a path of length less than \(L\) in
\(\Gamma\) from its smaller to its larger element. Each step of the
path is a reflection \(z\mapsto m-z\) with \(m\in D^+\) or a translation
by \(\pm h\) with \(h\in D^-\). After padding with identity maps to
length \(L-1\), there are at most \(M\) sequences of steps. The
composition of the steps is either a translation, by a positive amount
since the path goes from the smaller to the larger element, which
specifies the difference of the pair, or a reflection, which specifies
its sum. A Sidon set has at most one pair with a given positive
difference and at most one pair with a given sum, so each sequence of
steps arises from at most one short pair in \(A\). The estimate
\(M=n^{o(1)}\) follows from Lemma~\ref{lem:sidon-active}(a).
\end{proof}

\begin{lemma}\label{lem:sidon-reduction}
We have
\[
  |\mathcal A(U)|\le2^{o(\sqrt n)}\max_{C\in\mathcal C}
  \prod_{P\in\mathcal P(C)}(1+|P|).
\]
\end{lemma}

\begin{proof}
Since \(L\ge8r^2\) for large \(n\), Lemma~\ref{lem:sidon-short} shows
that Theorem~\ref{thm:deferred} applies with \(M=n^{o(1)}\).
Since \(|U|=O(\sqrt n)\) and \(L\to\infty\), the exponential factor
in~\eqref{eq:deferred-count} is \(2^{o(\sqrt n)}\).
Now use~\eqref{eq:sidon-product}.
\end{proof}

\subsection{A weighted Erd\H{o}s--Tur\'an bound}\label{sec:sidon-et}

It remains to bound the product in Lemma~\ref{lem:sidon-reduction}.
Erd\H{o}s and Tur\'an bounded the size of a Sidon set \(A\subseteq[n]\)
by smoothing its indicator function and using that every positive
difference occurs at most once in \(A\). We run their argument with the
weight \(|P|^{-1/2}\) on every element of a small part \(P\), so that the
part has total weight \(\sqrt{|P|}\). Elements of distinct small parts
are not adjacent in \(\Gamma\), so the quadruples they form are active.
By Lemma~\ref{lem:sidon-containers} there are few such quadruples, and
consequently the weighted number of pairs at a given difference is
essentially at most one. For \(C\in\mathcal C\), let
\[
  \Phi(C)=\sum_{P\in\mathcal P(C)}\sqrt{|P|}.
\]

\begin{lemma}\label{lem:weighted-et}
Every \(C\in\mathcal C\) satisfies \(\Phi(C)\le(1+o(1))\sqrt n\).
\end{lemma}

\begin{proof}
Give each \(x\) in a small part \(P\) the weight \(w(x)=|P|^{-1/2}\), and
let \(w(x)=0\) for all other integers \(x\), so that
\(\sum_xw(x)=\Phi(C)\). For \(h\ge1\), let
\(\rho(h)=\sum_xw(x)w(x+h)\). For a Sidon set with unit weights we would
have \(\rho(h)\le1\) for every \(h\). We first show that the total excess
over this bound is small, and then run the argument of Erd\H{o}s and
Tur\'an.

\begin{claim}\label{clm:sidon-excess}
The sum \(R\) of \(\rho(h)-1\) over all \(h\ge1\) with \(\rho(h)>1\)
satisfies \(R=o(n)\).
\end{claim}

\begin{proof}
Fix \(h\notin D^-\). Let \(q_h\) be the number of active quadruples in
\(C\) with difference \(h\), and let \(a_h\) be the number of three-term
progressions in \(C\) with common difference \(h\). We show that
\begin{equation}\label{eq:sidon-rho}
  \rho(h)\le1+q_h+a_h.
\end{equation}
Group the pairs \((x,x+h)\) with \(w(x)w(x+h)>0\) according to the
unordered pair \(\{P,P'\}\) of small parts containing \(x\) and \(x+h\).

Consider first a single group. If \(P=P'\), the group has
fewer than \(|P|\) pairs, each contributing \(1/|P|\), so it contributes
less than one. Suppose that \(P\neq P'\). Let \(S\subseteq P\) be the set
of left endpoints of the pairs going from \(P\) to \(P'\), and
\(S'\subseteq P\) the set of right endpoints of the pairs going from
\(P'\) to \(P\). The group has \(|S|+|S'|\) pairs, each contributing
\(1/\sqrt{|P||P'|}\). Every \(z\in S\cap S'\) is the middle term of the
progression \((z-h,z,z+h)\), whose two pairs belong to the group, and so
\(|S|+|S'|\le|P|+|S\cap S'|\). Arguing in the same way with the endpoints
in \(P'\), we see that the group contributes at most
\(\min\{|P|,|P'|\}/\sqrt{|P||P'|}\le1\) plus the number of progressions
with difference \(h\) formed by two of its pairs.

Next, take one pair \((x,x+h)\) from one group, and one pair
\((y,y+h)\) from each of the other groups. If the two pairs share an
element, then they form a three-term progression with common difference
\(h\). Otherwise they form a nondegenerate additive quadruple with
differences \(h\) and \(|x-y|\) and sum \(x+y+h\), and this quadruple is
active. Indeed, \(h\notin D^-\). If \(|x-y|\in D^-\), then \(\Gamma\)
joins \(x\) to \(y\) and \(x+h\) to \(y+h\); if \(x+y+h\in D^+\), then
\(\Gamma\) joins \(x\) to \(y+h\) and \(x+h\) to \(y\). By~(P2), in
both cases the two pairs would lie in the same group.

The configurations found in the last two paragraphs are distinct, since
each of them determines its two pairs, and they are progressions or
active quadruples with difference \(h\). If there are \(g\ge1\) groups,
we have found \(g-1\) such configurations between distinct groups, and
the bound on the contribution of each group gives~\eqref{eq:sidon-rho}.

For \(h\in D^-\) we use
\(\rho(h)\le\sum_xw(x)^2=|\mathcal P(C)|\le|C|\). Each quadruple has two
differences, and a progression is determined by its middle term \(a\)
and a pair with sum \(2a\), so
\[
  R\le|D^-|\,|C|+2\Qact(C)+\sum_{a\in C}r_C^+(2a).
\]
The three terms are \(O(\sqrt n(\log n)^2)\) by
Lemma~\ref{lem:sidon-active}(a), \(o(n)\) by~\eqref{eq:active-few}, and
\(O(n^{3/4})\) by~\eqref{eq:progressions}.
\end{proof}

To finish, we apply the smoothing argument of Erd\H{o}s and Tur\'an. For
an integer \(\ell\ge1\), the function \(g(y)=\sum_{i=0}^{\ell-1}w(y-i)\)
has sum \(\ell\Phi(C)\) and support of size at most \(n+\ell\). By the
definition of \(R\), we have
\(\sum_{h=1}^{\ell-1}(\ell-h)\rho(h)\le\binom\ell2+\ell R\), and hence
\[
  \sum_yg(y)^2=\ell\sum_xw(x)^2+2\sum_{h=1}^{\ell-1}(\ell-h)\rho(h)
  \le\ell|C|+\ell^2+2\ell R.
\]
By the Cauchy--Schwarz inequality,
\[
  \Phi(C)^2\le(n+\ell)\Bigl(1+\frac{|C|+2R}{\ell}\Bigr).
\]
Taking \(\ell=\lceil\sqrt{n(|C|+2R+1)}\rceil\), which is \(o(n)\), gives
\(\Phi(C)^2\le(1+o(1))n\).
\end{proof}

\begin{proof}[Proof of Proposition~\ref{prop:terminal}]
We first show that
\begin{equation}\label{eq:log-sqrt}
  \log(1+a)\le\tfrac12\log5\cdot\sqrt a
  \qquad\text{for every positive integer }a,
\end{equation}
that is, \((1+a)^2\le5^{\sqrt a}\). For \(a\le5\) this is checked
directly, with equality at \(a=4\). For \(x\ge\sqrt5\), the derivative of
\(x\log5-2\log(1+x^2)\) is
\[
  \log5-\frac{4x}{1+x^2}\log e\ge\log5-\frac{2\sqrt5}3\cdot\frac32
  =\log5-\sqrt5>0,
\]
since \(4x/(1+x^2)\) is decreasing for \(x\ge1\) and \(\log e<3/2\). This
proves~\eqref{eq:log-sqrt} for \(a\ge5\).

For every \(C\in\mathcal C\), \eqref{eq:log-sqrt} and
Lemma~\ref{lem:weighted-et} give
\[
  \sum_{P\in\mathcal P(C)}\log(1+|P|)
  \le\tfrac12\log5\cdot\Phi(C)
  \le\bigl(\tfrac12\log5+o(1)\bigr)\sqrt n,
\]
and Lemma~\ref{lem:sidon-reduction} proves the proposition. The
maximizer \(a=4\) mirrors the lower-bound construction: there too
each element of the modular Sidon set has four representatives, giving
\(1+4=5\) choices.
\end{proof}

\section{Further applications}\label{sec:further}

Deferral also gives the three results below. Their proofs are complete,
but the manuscripts containing them have not yet been revised for
publication. They are available on the author's webpage,
\url{https://tianchi-yang.github.io/}, and will be updated there.

\subsection{Induced Ramsey numbers}

For graphs \(F\) and \(H\), the \emph{induced Ramsey number}
\(\rind(F,H)\) is the least \(n\) for which there is a graph on \(n\)
vertices such that every red--blue colouring of its edges contains an
induced copy of \(F\) whose edges are all blue or an induced copy of \(H\)
whose edges are all red. Erd\H{o}s~\cite{erdos84} conjectured that
\(\rind(H,H)\le2^{O(t)}\) for every graph \(H\) on \(t\) vertices, and
Arag\~ao, Campos, Dahia, Filipe and Marciano~\cite{acdfm} recently proved
this conjecture using a conditional form of the container lemma of Campos
and Samotij~\cite{cs}. For graphs \(F\) and \(H\) on \(k\) and \(t\)
vertices, Kohayakawa, Pr\"omel and R\"odl~\cite{kpr} conjectured that
\(\rind(F,H)\) is polynomial in \(t\) for every fixed \(F\). Zhong, Kashima
and Mao~\cite{zkm} and, independently, Andrades, Campos, Filipe and
Kats~\cite{and} proved this in the form \(\rind(F,H)\le t^{O(k^2)}\);
Andrades, Campos, Filipe and Kats also asked whether the exponent can be
made linear in \(k\). Combining deferral with a signed form of the
container lemma of~\cite{acdfm}, we prove the following.

\begin{theorem}\label{thm:ramsey}
For all graphs \(F\) and \(H\) with \(2\le k=|V(F)|\le t=|V(H)|\),
\[
  \rind(F,H)\le(2t/k)^{10^8k}.
\]
\end{theorem}

Up to the constant in the exponent, this is best possible even for
complete graphs: \(\rind(K_k,K_t)\) is the Ramsey number \(R(k,t)\), and
the first-moment criterion of Spencer~\cite{spencer} gives
\(\log R(k,t)\ge\frac1{16}k\log(2t/k)\) for \(2\le k\le t\).

\subsection{Sets without arithmetic progressions}

Let \(a_k(n)\) be the number of subsets of \([n]\) without a nonconstant
\(k\)-term arithmetic progression, and let \(r_k(n)\) be the maximum size
of such a subset. Balogh, Liu and Sharifzadeh~\cite{bls} proved that
\(\log a_k(n)=O_k(r_k(n))\) for infinitely many \(n\), and Morris,
Ortega and Ru\'e~\cite{mor} proved this bound for all \(n\) and obtained
the constant one for infinitely many \(n\) when \(k\ge5\).

\begin{theorem}\label{thm:ap}
For every fixed \(k\ge3\),
\[
  a_k(n)=2^{(1+o(1))r_k(n)}
\]
as \(n\to\infty\) through a set of logarithmic density one that can be
chosen independently of \(k\).
\end{theorem}

\subsection{List colourings of linear hypergraphs}

A hypergraph is \emph{linear} if any two distinct edges share at most one
vertex. The list chromatic number \(\chi_{\mathrm{list}}(\mathcal G)\) of
a hypergraph \(\mathcal G\) is the least \(m\) such that, for every
assignment of lists of \(m\) colours to the vertices, there is a proper
colouring in which every vertex receives a colour from its list. Saxton
and Thomason~\cite{st} proved that every linear \(r\)-uniform hypergraph
\(\mathcal G\) of average degree \(d\) satisfies
\(\chi_{\mathrm{list}}(\mathcal G)\ge(\frac1{(r-1)^2}+o(1))\log_rd\), and
that \(\chi_{\mathrm{list}}(\mathcal G)\ge(\frac1{r-1}-o(1))\log_rd\) if
\(\mathcal G\) is in addition \(d\)-regular; they asked whether the second
bound holds with average degree in place of regularity.

\begin{theorem}\label{thm:list}
Let \(r\ge3\) be fixed. Every linear \(r\)-uniform hypergraph
\(\mathcal G\) of average degree \(d\) satisfies
\[
  \chi_{\mathrm{list}}(\mathcal G)\ge\Bigl(\frac1{r-1}-o(1)\Bigr)\log_rd,
\]
where \(o(1)\to0\) as \(d\to\infty\), uniformly over all such
hypergraphs.
\end{theorem}

\section*{AI disclosure}

In earlier work on \(C_4\)-free graphs~\cite{yang}, the author deferred
constraints from heavily overlapping configurations until the final count
within each container. This paper develops that idea into a general method.

The author completed the proofs with assistance from GPT~6.1 and
subsequently used GPT, Claude, and Gemini to polish the exposition.
The author checked all suggestions made by these tools and takes full
responsibility for the mathematical content of the paper.

\end{document}